\documentclass{article}

\usepackage{amsmath}
\usepackage[all]{xy}

\usepackage{graphicx}
\usepackage{amssymb,amsthm}

\newtheorem{theorem}{Theorem}[section]
\newtheorem{corollary}[theorem]{Corollary}
\newtheorem{example}[theorem]{Example}

\newtheorem{lemma}[theorem]{Lemma}
\newtheorem{proposition}[theorem]{Proposition}
\newtheorem{remark}[theorem]{Remark}

\newtheorem{definition}[theorem]{Definition}

\newcommand{\ie}{i.e.\ }
\newcommand{\eg}{e.g.\ }

\newcommand{\defn}{\emph}
\newcommand{\Z}{\mathbb{Z}}

\newcommand{\R}{\mathbb{R}}
\newcommand{\C}{\mathbb{C}}

\newcommand{\PP}{\mathbb{P}}

\newcommand{\cat}[1]{\mathbf{#1}}

\newcommand{\map}[2]{{\rm Map}\left(#1,#2\right)}
\newcommand{\pomap}[2]{{\rm Map}^{\leq}\!\left(#1,#2\right)}

\newcommand{\holink}[2]{{\rm holink}\left(#1,#2\right)}

\newcommand{\funct}[2]{{\left[#1,#2\right]}}

\newcommand{\open}[1]{\cat{U}(#1)}

\newcommand{\et}[1]{\cat{Et}\!\downarrow\!#1}
\newcommand{\br}[1]{\cat{Br}\!\downarrow\!#1}

\newcommand{\shc}[1]{\cat{Sh^c}\!\downarrow\!#1}
\newcommand{\cshc}[1]{\cat{Cosh^c}\!\downarrow\!#1}

\newcommand{\ulpo}[1]{\cat{UL^{po}}\!\downarrow\!#1}
\newcommand{\ulop}[1]{\cat{UL^{op}}\!\downarrow\!#1}

\renewcommand{\cosh}[1]{\mathcal{#1}}

\begin{document}

\title{The fundamental category of a stratified space}          

\author{Jon Woolf}       
\date{September 2013}      

\maketitle

\section{Introduction}

Covers of a (nice) topological space $X$ are classified by the fundamental group $\pi_1X$ or, if we wish to avoid assuming that $X$ is connected and has a basepoint, by the fundamental groupoid $\Pi_1X$. If $X$ is a stratified space then it is natural to allow the covers to be stratified too. MacPherson observed (unpublished) that local homeomorphisms onto $X$ which are covers when restricted to each stratum are classified by a modified version of the fundamental groupoid. The objects of this are the points of $X$ and the morphisms are homotopy classes of paths which `wind outwards from the deeper strata', \ie once they leave a stratum they do not re-enter it. The notion of homotopy here is that of  a family of such paths. MacPherson phrased his result in terms of constructible sheaves rather than covering spaces, but the notions are equivalent under the well-known correspondence of sheaves and their \'etale spaces.

MacPherson's ideas were published and extended by Treumann \cite{treumann} in which he refers to paths which wind outwards as exit paths. For any topologically stratified space he defines the `exit $1$-category' --- objects are points, morphisms are homotopy classes of exit paths (with a tameness assumption on the homotopy) --- and the `exit $2$-category' --- objects are points, $1$-morphisms are exit paths and $2$-morphisms `tame' homotopies through exit paths. He shows that the set-valued functors on the exit $1$-category are equivalent to constructible sheaves and that category-valued functors on the exit $2$-category  are equivalent to constructible stacks. Interesting examples of constructible stacks are given by categories of perverse sheaves on $X$. However, there does not yet seem to be any way of identifying the particular representations which correspond to these. 

This paper develops these ideas by defining a `fundamental category' for any `pre-ordered space'. Stratified spaces are particular examples, and in this case our definition reduces to Treumann's exit $1$-category (but without the tameness condition on homotopies). The fundamental category has good properties for a very wide class of stratified spaces, namely homotopically stratified sets with locally  $0$ and $1$-connected strata. For these:
\begin{itemize}
\item The fundamental category can be computed in terms of homotopy groups of the strata and of the homotopy links --- see \S \ref{computing}.
\item The fundamental groupoid can be recovered by localising the fundamental category at all morphisms --- see Corollary \ref{representations}.
\item  \label{gen Mac} Covariant set-valued functors from the fundamental category classify constructible sheaves and, dually, contravariant functors classify constructible cosheaves --- see Theorem \ref{fund groupoid}. Geometrically these can be interpreted respectively as classifications of stratified \'etale and branched covers.
\end{itemize}
This class of spaces includes Whitney stratified spaces, Thom-Mather stratified spaces, topologically stratified spaces and Siebenmann's locally cone-like spaces. In particular we recover MacPherson's result. We do not consider the analogue of Treumann's exit $2$-category but it seems likely that his results would generalise to homotopically stratified sets. 

\subsection{Structure of the paper}

A pre-ordered space, or po-space, has a distinguished subset of paths: the po-paths are those paths which are also order-preserving maps, where $[0,1]$ is ordered by $\leq$. When $X$ is a stratified space po-paths are precisely Treumann's exit paths. The fundamental category $\Pi^{po}_1X$ of a po-space $X$, defined in \S\ref{fund-cat}, is an ordered analogue of the fundamental groupoid in which morphisms are given by homotopy classes of po-paths. 

In \S\ref{hss} we restrict the discussion  to a homotopically nice class of filtered spaces, Quinn's homotopically stratified sets.  These spaces have two advantages for our purposes. Firstly they are very general, subsuming almost any other notion of stratified space. Secondly, the language used in their definition, particularly that of the homotopy link, is well-suited for talking about po-paths. In particular it allows us to reduce to the simpler case of a space with only two strata. In \S\ref{computing} we explain an approach to computing the fundamental category based on this reduction.

The fundamental groupoid of a locally $0$ and $1$-connected space classifies covers, equivalently local systems. In \S\ref{stratified covers} we show that the fundamental category plays a similar r\^{o}le for homotopically stratified spaces with locally $0$ and $1$-connected strata. However, because po-paths are not necessarily reversible, we obtain two different classification results. Covariant set-valued representations of the fundamental category of a homotopically stratified set correspond to constructible sheaves, or equivalently to stratified \'etale covers. Dually, contravariant representations correspond to constructible cosheaves, or equivalently to stratified branched covers (which we define in terms of Fox's notion of a complete spread \cite{fox}). 

As  a corollary of this classification result we deduce that the fundamental groupoid of a homotopically stratified set with locally $0$ and $1$-connected strata can be recovered by localising the fundamental category at the set of all morphisms: $\Pi_1X$ is the  `groupoidification' of  $\Pi^{po}_1X$.

In \S\ref{config} we consider an example, the symmetric product $SP^n\C$ with the natural stratification indexed by partitions of $n$. Morphisms in the fundamental category can be expressed in terms of various subgroups of the braid group $B_n$ and symmetric group $S_n$ associated to partitions. 

Appendix \ref{pre-orders} explains how filtered spaces arise as po-spaces with a certain natural compatibility between the pre-order and the topology. This is included to make the case that filtered and stratified spaces are not exotic examples in the world of ordered topology but the bread and butter of the subject. Finally, since cosheaves and complete spreads are less familiar than sheaves and \'etale maps we give a brief review of the relevant theory in Appendix \ref{cosheaves}.


\subsection{Related work}

We have already discussed the relation to MacPherson and Treumann's ideas. In a less geometric vein, there has been work on homotopy theory for ordered or directed spaces in category theory and theoretical computer science. Unfortunately there is a plethora of slightly different definitions. In \cite{Kahl:kl} Kahl shows that the category of spaces equipped with partial orders and maps between them has a closed model structure with {\em unordered} homotopies as the notion of homotopy. Bubenik and Worytkiewicz \cite{Bubenik:so} follow a similar programme for {\em locally ordered} spaces --- spaces $X$ with a partial order which need only be transitive locally, but which is closed as a subset of $X \times X$. However, here it is only possible to show that  such spaces are {\em contained} in a closed model category. Finally, Grandis \cite{M.-Grandis:2003qz} considers a more general notion of {\em directed} spaces --- spaces equipped with a suitable class of `directed' paths. He shows that the directed category has good properties (existence of limits and colimits, exponentiable directed interval etc) and develops directed homotopy theory within it using {\em ordered} homotopies (so that directed homotopy is not an equivalence relation). The intended application in all three cases is in theoretical computer science, to the theory of concurrent systems. Grandis also uses directed homotopy to produce interesting examples of higher-dimensional categories, see \cite{M.-Grandis-8-2006:2006dw}. 

\subsection{Acknowledgments}

This work was funded by the Leverhulme Trust (grant reference: F/00 025/AH). I am grateful to Michael Loenne who explained the ideas behind the example in \S\ref{config} and to Beverley O'Neill (generously funded by a Nuffield Undergraduate Bursary) who worked out the details. 
I am indebted to David Miller for sparing my blushes by pointing out a serious error in an earlier preprint. Thanks also to Andrei Prasolov for bringing my attention to a mistake in Appendix B in the originally published version of this paper, where too strong a result relating cosheaves to spreads was claimed and erroneously proved. I would like to thank Ivan Smith and particularly Tom Leinster for helpful discussions and suggestions. And finally, I would like to thank the referee for his meticulous reading and helpful corrections and comments.

\section{The fundamental category}
\label{fund-cat}

There are various notions of `ordered space'. The one we use, the notion of a po-space, is the simplest. It is a topological space with a pre-order on the set of points. Recall that a \defn{pre-order} is a set equipped with a reflexive and transitive relation $\leq$ and that a map is \defn{increasing} if $p \leq q \Rightarrow f(p) \leq f(q)$. Po-spaces form a category with maps between them being both continuous and increasing. We will refer to these as po-maps for brevity.

Notice that there need be no compatibility between the topology and the pre-order. However, if we impose a natural compatibility condition then the resulting po-spaces are filtered spaces (see Appendix \ref{filtered}). These are the  po-spaces which arise most often in geometry and topology; all the examples we consider will be of this kind, indeed they will be stratified spaces (see \S\ref{hss}). The relation between po-spaces and filtered spaces is explained more fully in Appendix \ref{alexandrov}. 

From now on we will assume that all spaces are compactly generated.  Spaces of maps are topologised with the $k$-ification of the compact-open topology (so that they too are compactly generated).

For po-spaces $X$ and $Y$ let $\pomap{X}{Y}$ be the set of po-maps between them. We topologise this as a subspace of  $\map{X}{Y}$. Let $I$ be the ordered interval, \ie $[0,1]$ equipped with the standard order. An element of $\pomap{I}{X}$ is a continuous path $\gamma:[0,1]\to X$ such that $\gamma(s) \leq \gamma(t)$ whenever $s\leq t$; we call it a \defn{po-path} in $X$. The start and end of a po-path determine a continuous map
$$
E_0\times E_1 : \pomap{I}{X} \to X^2 : \gamma \mapsto \left(\gamma(0) , \gamma(1) \right).
$$

\begin{definition}
The fundamental category $\Pi^{po}_1 X$ of  a po-space $X$ is the category whose objects are the points of $X$ and whose morphisms from $x$ to $x'$ are the (path) connected components $\pi_0 (E_0\times E_1)^{-1}(x,x')$ of the space of po-paths from $x$ to $x'$. That is, a morphism from $x$ to $x'$ is a homotopy class of po-paths from $x$ to $x'$ where the homotopy is through po-paths. Composition is defined by concatenation of po-paths. The fundamental category is functorial: a po-map $f: X \to Y$ induces a functor $\Pi^{po}_1(f) : \Pi^{po}_1X \to \Pi^{po}_1Y$. 
\end{definition}

\begin{example}
Let $P$ be a poset with the descending chain condition, \ie any descending chain of elements of $P$ is eventually constant. Consider $P$ as a po-space by giving it the Alexandrov topology (see Appendix \ref{pre-orders}). Any po-path in $P$ is of the form
$$
\gamma(t) = p_i \ \textrm{for}\ t \in (t_{i-1}, t_{i}] \cap [0,1] \qquad i=0,\ldots,n
$$
where $p_0 \leq \cdots \leq p_n$ is a chain in $P$ and $t_{-1} <0 < t_0 < \cdots < t_n = 1$ an increasing sequence in $\R$. There is a homotopy through po-paths
$$
\eta(s,t) = 
\left\{
\begin{array}{ll}
\gamma(t) & s=0\\
p_0 & t=0\\
p_n & \textrm{otherwise},
\end{array}
\right.
$$
from $\gamma$ to the po-path starting at $p_0$ and moving instantly to $p_n$. Hence $\Pi^{po}_1 P \cong P$ thought of as a category with a single morphism from $p$ to $p'$ when $p \leq p'$.

\end{example}

\section{Homotopically stratified sets}
\label{hss}

A stratified space is a filtered space together with some information on how the strata glue together. There are many ways in which we can specify glueing data and hence many types of stratified space. In a homotopy-theoretic context the most flexible is Quinn's notion of a homotopically stratified set (defined below). This is the notion of stratified space we will use. It is very general, in particular Siebenmann's locally cone-like stratified spaces, Thom-Mather stratified spaces and Whitney stratified spaces are all homotopically stratified sets. The main example in this paper is the symmetric product $SP^n\C$. We stratify this with one stratum for each partition of $n$; the corresponding stratum is the subset of configurations in which the $n$ points coalesce according to the partition (see \S \ref{config}). 

The fundamental category of a homotopically stratified set is equivalent to a similar category in which morphisms are homotopy classes of po-paths which only pass through one or two strata. This allows us to describe morphisms in terms of the homotopy groups of strata and of homotopy links of pairs of strata.

In order to give the definition of a homotopically stratified set we introduce some terminology. Suppose $B$ is a filtered space. A subspace $A \subset B$ is \defn{tame} if it is a nearly stratum-preserving deformation retract of a neighbourhood $N$ of $A$, \ie there is a deformation retraction of $N$ onto $A$ such that points remain in the same stratum under the retraction until the last possible moment, when they must flow into $A$. The \defn{homotopy link} $\holink{B}{A}$ of a subset $A$ of $B$ is the space of paths (equipped with the compact-open topology) in $B$ with $\gamma(0) \in A$ and $\gamma(0,1] \subset B - A $, \ie the space of paths which start in $A$ but leave it immediately. Evaluation at $t$ defines a map $E_t : \holink{A \cup B}{A} \to A\cup B$. 
\begin{definition}
A \defn{homotopically stratified set} is a filtered space $X$ with finitely many connected strata $X^i$ such that for any pair $i\leq j$
\begin{enumerate}
\item  \label{tame subspace} the inclusion $X^i \hookrightarrow X^i\cup X^j$ is tame and
\item \label{fibred holinks} the evaluation map  $E_0 : \holink{X^i\cup X^j}{X^i} \to X^i$ is a fibration.
\end{enumerate}
\end{definition}
It is not immediately apparent why this is a good definition. One reason is that, if we assume that $X$ is a metric space, then 
$$
\xymatrix{
\holink{X^i\cup X^j}{X^i} \ar[d]_{E_0} \ar[r]^{\qquad \quad E_1} & X^j \ar@{^{(}->}[d]\\
X^i \ar@{^{(}->}[r] & X^i \cup X^j
}
$$ 
is a homotopy push-out \cite[Lemma 2.4]{quinn}. Intuitively the homotopy link plays the r\^{o}le of the boundary of a regular neighbourhood of $X^i$ in $X^i\cup X^j$.

A homotopically stratified set is naturally a po-space, where we give it the pre-order corresponding to the underlying filtration. Homotopy links can then be described in terms of po-maps. Let $I_1$ be the interval $[0,1]$ made into a po-space via the filtration $\{0\} \subset [0,1]$ --- see Appendix \ref{alexandrov}.  Then
$$
\holink{X^i\cup X^j}{X^i} = (E_0\times E_1)^{-1} \left(X^i \times X^j\right) \subset \pomap{I_1}{X}.
$$ 
In fact it is a good heuristic when working with homotopically stratified sets that their homotopy theory can be understood in terms of maps from $I_1$. The next lemma, which shows that po-paths are homotopic to elements of the holink, is an illustration of this principle. 
\begin{lemma}
Let $X$ be a homotopically stratified set and $\gamma:[0,1]\to X$ a po-path in $X$ from $x_i\in X^i$ to $x_j\in X^j$. Then there is a homotopy of $\gamma$, relative to its end points, to a po-path $\tilde{\gamma} \in \holink{X^i\cup X^j}{X^i}$. Moreover, $\tilde{\gamma}$ is unique up to homotopy through po-paths in the homotopy link. 
\end{lemma}
\begin{proof}
In general the po-path $\gamma:[0,t] \to X$ will pass through several strata. Let $t_{-1}< 0\leq t_0 <  \cdots < t_{n-1} < t_n=1$ be such that $\gamma(t)$ is in one stratum for $t\in (t_{k-1},t_k] \cap [0,1]$ for $k=0,1,\ldots,n$. In particular the trace of $\gamma$ is in $X^i$ for $t\in[0,t_1]$ and in $X^j$ for $t\in (t_n,1]$. There is a homotopy of $\gamma$ which moves it off the intermediate strata one by one starting from the highest. The segment $\gamma : (t_{n-1},t_n] \to X$ is a path in some stratum $X^k$. The final segment $\gamma : [t_n,1] \to X$ is a lift (inverse image) of $\gamma(t_n)$ along the start point map $E_0 : \holink{X^k\cup X^j}{X^k} \to X^k$. By definition $E_0$ is a fibration so we can extend this lift along $\gamma: (t_{n-1},t_n] \to X^k$. The result is a map 
$$
\eta : (t_{n-1},t_n] \times [t_n,1] \to X
$$
with $\eta(-,t_n) = \gamma(-)$, $\eta(t_n,-) = \gamma(-)$ and $\eta(s,t) \in X^j$ for $t\neq t_n$. This provides a homotopy, through po-paths, between $\gamma$ and 
$$
\gamma_1(t) =\left\{
\begin{array}{ll}
 \gamma(t) &  0\leq t \leq t_{n-1}   \\
 \eta(t_{n-1},t_n-t_{n-1}+t) &  t_{n-1}< t \leq t_{n-1} + 1-t_n   \\
 \eta(t_n-1+t,1) &  t_{n-1} + 1-t_n< t\leq 1.   
\end{array}
\right.
$$
The po-path $\gamma_1$ now passes through one fewer intermediate strata (since it avoids $X^k$). Continuing inductively we obtain the desired $\tilde{\gamma}$. Uniqueness up to homotopy follows from multiple applications of the uniqueness up to homotopy of the extension of a lift along a fibration. The situation is most easily apprehended by looking at the following diagram.
\begin{center}
   \includegraphics[width=4.5in]{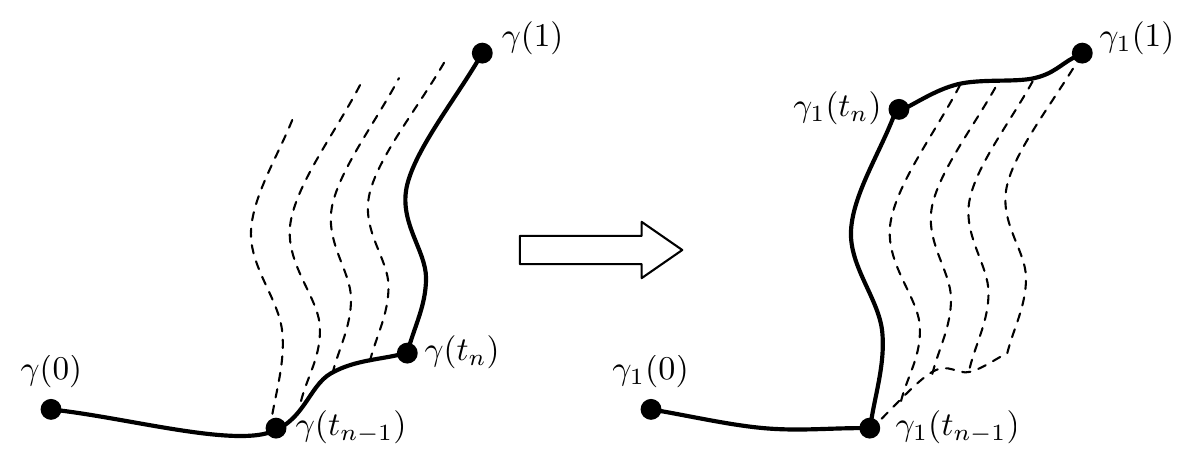} 
   \end{center}
\end{proof}

Of course we would really like a relative version of this result, allowing us to find homotopies of families of po-paths to families in the holink. This does not follow for free from the above lemma; the main difficulty is that the point at which po-paths in a family leave a given stratum is not a continuous function of the parameters of the family. Nevertheless Miller has recently announced
\begin{theorem}[See {\cite[Theorem 3.9]{miller}}.]
\label{miller}
Let $X$ be a homotopically stratified set and suppose $X$ is a metric space. Then the space of po-paths beginning in $X^i$ and ending in $X^j$ is homotopy equivalent to $\holink{X^i\cup X^j}{X^i}$ by a homotopy which fixes start and end points and takes po-paths `instantly' into the homotopy link.
\end{theorem}
This theorem allows us to give an equivalent, but simpler, definition of the fundamental category. Let $\Pi^{ho}_1X$ be the category with objects the points of $X$ and with morphisms  from $x\in X^i$ to $x'\in X^j$ given by the (path) connected components of
$$
(E_0 \times E_1)^{-1}(x,x') \subset \holink{X^i\cup X^j}{X^i}.
$$
In order to define composition we need to concatenate paths and then choose a homotopic (through po-paths) path in the holink. The above theorem guarantees the result is well-defined, up to homotopy through po-paths.
\begin{corollary}
The functor $\Pi^{ho}_1X \to \Pi^{po}_1X$ induced by the inclusion of the homotopy link in the space of po-paths is an equivalence.
\end{corollary}
The proof is immediate from Theorem \ref{miller}. Later, in Corollary \ref{ho to po}, we will give an independent proof of this in the case when the strata of $X$ are locally $0$ and $1$-connected. This result allows us to reduce to the two-stratum case: morphisms in the fundamental category of a homotopically stratified set can be described in terms of the homotopy links, they do not depend on any intermediate strata.

\subsection{Computing the fundamental category}
\label{computing} 
In this section we give several fibration sequences for computing morphisms in the fundamental category of a homotopically stratified set $X$.
 \begin{lemma}
\label{fibration lemma}
Write $H^{ij}$ for the homotopy link $\holink{X^i\cup X^j}{X^i}$. Fix a basepoint $x_i$ for each stratum $X^i$. The four downward maps in the diagram
$$
\xymatrix{
E_1^{-1}(x_j) \ar@{^{(}->}[r] \ar[d]_{E_0} &H^{ij} \ar[dl]^{E_0} \ar[dr]_{E_1} & E_0^{-1}(x_i) \ar@{_{(}->}[l] \ar[d]^{E_1}\\
X^i && X^j
}
$$
are fibrations. In addition $E_0 \times E_1 : H^{ij}  \to X^i \times X^j$ is a fibration. 
\end{lemma}
\begin{proof}
Since $X$ is homotopically stratified $E_0 : H^{ij}  \to X^i$ is a fibration. Now consider the restriction $E_1^{-1}(x') \to X^i$ of $E_0$ to paths ending at $x'$. Given a family $F: A \times [0,1] \to X^i$ of paths in $X^i$ and a lift $G: A \times \{0\} \to E_1^{-1}(x')$ of the start points we have a lift $G: A \times [0,1] \to H^{ij}$. Then the `diagonal' family $\widetilde{F} : A \times [0,1] \to E_1^{-1}(x')$ given by
$$
\widetilde{F}(a,s)(t) = G\left(a, s(1-t) \right) (t)
$$
is a lift of $F$ to $E_1^{-1}(x')$ along $E_0$ starting at $G(-,0)$. Hence the restriction of $E_0$ to $E_1^{-1}(x') $ is a fibration.

The end point map $E_1 :H^{ij}  \to X^j$ is a fibration because, given a family $F:A \times [0,1] \to X^j$ of paths in $X^j$ and a lift $G: A \times \{0\} \to H^{ij}$ of the starting points, the family given by the composition of paths
$$
\widetilde{F}(a,t) =  G(a)(-) \circ F(a,-)|_{[0,t]} 
$$
is a lift of $F$ to $H^{ij}$. It is easy to see that this construction of a lift also shows that the restriction of $E_1$ to $E_0^{-1}(x_i)$ is a fibration.

Finally, we must show that $E_0\times E_1 :H^{ij} \to X^i \times X^j$ is a fibration. If $F = (F_1,F_2): A \times [0,1] \to X^i \times X^j$ is a family of paths and $G : A \times \{0\} \to H^{ij}$ a lift of the starting points then, using essentially the same argument which showed that $E_0:E_1^{-1}(x')\to X^i$ is a fibration, we can find a lift $G: A\times [0,1] \to H^{ij}$ of $F_1$ along $E_0$ such that $G(a,s)(1) = G(a,0)(1) = F_2(a,0)$. Then the composition of paths
$$
\widetilde{F}(a,s) = G(a,s) \circ F_2(a,-)|_{[0,s]}
$$
is a lift of $F$ along $E_0\times E_1$.
\end{proof}

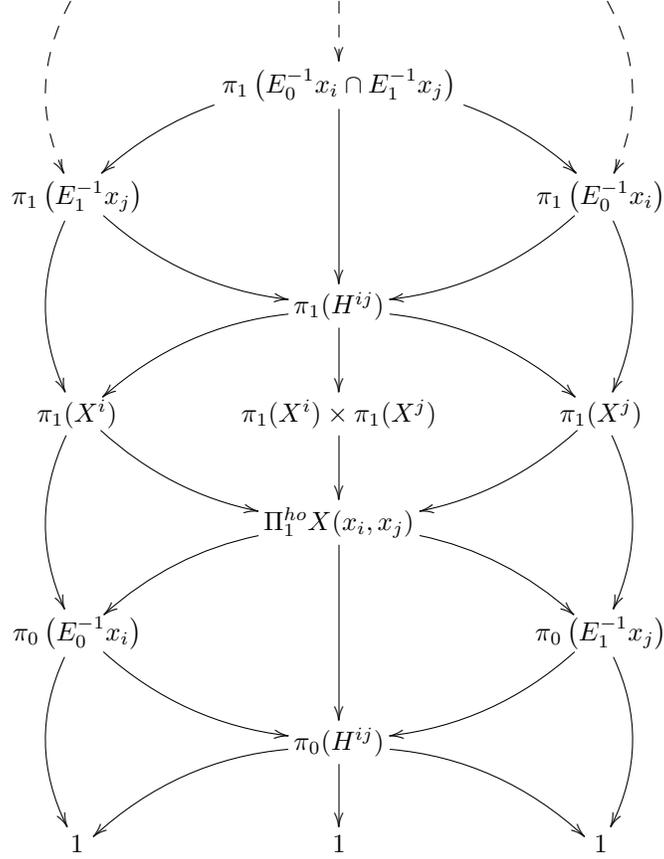
\begin{figure}[htbp]
\begin{center}
$$
\xymatrix{
 \ar@/_1pc/@{-->}[dd] &\ar@{-->}[d] & \ar@/^1pc/@{-->}[dd]  \\
& \pi_1\left( E_0^{-1}x_i \cap E_1^{-1}x_j  \right) \ar@/^1pc/[dr] \ar[dd] \ar@/_1pc/[dl]  & \\
\pi_1\left( E_1^{-1}x_j  \right) \ar@/_1pc/@<.5ex>[dr] \ar@/_1pc/[dd] &&
 \pi_1\left( E_0^{-1}x_i \right) \ar@/^1pc/@<-0.5ex>[dl] \ar@/^1pc/[dd] \\
& \pi_1(H^{ij}) \ar@/^1pc/@<-0.5ex>[dr] \ar[d] \ar@/_1pc/@<.5ex>[dl] &\\
\pi_1(X^i)  \ar@/_1pc/[dd] \ar@/_1pc/@<.5ex>[dr] & 
\pi_1(X^i) \times \pi_1(X^j) \ar[d] & 
 \pi_1(X^j)  \ar@/^1pc/[dd] \ar@/^1pc/@<-0.5ex>[dl] \\
 & \Pi^{ho}_1X(x_i,x_j) \ar@/^1pc/@<-0.5ex>[dr] \ar[dd] \ar@/_1pc/@<.5ex>[dl] & \\
\pi_0\left( E_0^{-1}x_i \right) \ar@/_1pc/@<.5ex>[dr] \ar@/_1pc/[dd] &&
 \pi_0\left( E_1^{-1}x_j   \right) \ar@/^1pc/@<-0.5ex>[dl] \ar@/^1pc/[dd] \\
& \pi_0(H^{ij}) \ar@/^1pc/@<-0.5ex>[dr] \ar[d] \ar@/_1pc/@<.5ex>[dl] &\\
1 &1 & 1
}
$$
\caption{A commutative braided diagram showing the (lower parts of the) long exact sequences arising from the fibrations in Lemma \ref{fibration lemma}. To aid reading we have suppressed all basepoints; they are the points $x_i\in X^i, x_j\in X^j$ or the path $\gamma^{ij}\in \holink{X^i\cup X^j}{X^i}$ from $x_i$ to $x_j$ as appropriate. Recall that $\Pi^{ho}_1X(x_i,x_j) \cong \pi_0\left( E_0^{-1}x_i \cap E_1^{-1}x_j \right)$ and that we assume strata are connected, hence the row of $1$s at the bottom.  See \S \ref{config} for an example of a computation using these sequences.
}
\label{braided LES}
\end{center}
\end{figure}

A compatible choice of basepoints for the spaces appearing in Lemma \ref{fibration lemma} is given by choosing a basepoint $\gamma^{ij} \in H^{ij}$ with $E_0(\gamma^{ij})=x_i$ and $E_1(\gamma^{ij})=x_j$ for each $i<j$, \ie a path with $\gamma^{ij}(0) =x_i, \gamma^{ij}(1)=x_j$ and $\gamma^{ij}(0,1] \subset X^j$. With this choice we obtain five long exact sequences corresponding to the five fibrations. These are displayed in a single commutative `braided' diagram in Figure \ref{braided LES}. 

The long exact sequences in Figure \ref{braided LES} give us several methods for computing sets of morphisms in $\Pi^{ho}_1X$. We would also like to be able to compute compositions of morphisms. In general this is fiddly; the complications arise in keeping track of the combinatorics when the homotopy links and their fibres are disconnected. Here we will treat only the simpler case when they are connected. In \S \ref{config} we treat a more difficult example where this fails. 

\begin{lemma} 
\label{composing}
Let $X$ be a homotopically stratified set. Assume that the fibres $E_0^{-1}(x_i)\subset H^{ij}$ are connected for each $i \leq j$ so that we have an exact sequence
$$
\xymatrix{
\cdots \ar[r] &  \pi_1E_0^{-1}(x_i) \ar[r]^{{\pi_1E_1}}&  \pi_1(X^j) \ar[r]^{\Gamma^{ij}}& \Pi^{ho}_1X(x_i,x_j) \ar[r] &  1
}
$$
where $\Gamma^{ij}$ is the map induced by pre-composition by $\gamma^{ij}$.  (For the sake of readability we have suppressed the basepoints.) Then for $i \leq j \leq k$
$$
\xymatrix{
\pi_1(H^{jk}) \times \pi_1(X^k) \ar[r]^{\pi_1E_1\times 1} \ar@{->>}[d]_{\pi_1E_0 \times 1}& \pi_1(X^k) \times \pi_1(X^k) \ar[d]^{{\rm compose}}\\
\pi_1(X^j) \times \pi_1(X^k)  \ar@{->>}[d]_{\Gamma^{ij} \times\Gamma^{jk}} & \pi_1(X^k) \ar@{->>}[d]^{\Gamma^{ik}} \\
\Pi_1^{ho}X(x_i,x_j) \times \Pi_1^{ho}X(x_j,x_k) \ar[r] & \Pi_1^{ho}X(x_i,x_k) 
}
$$
commutes, where the bottom map is composition in $\Pi^{ho}_1X$ and the bottom right map $\Gamma^{ik}$ is pre-composition  by $\gamma^{ij}\gamma^{jk}$.  Hence, at least in principle, we can compose a pair of morphisms by choosing a lift  to $\pi_1(H^{jk}) \times \pi_1(X^k)$ along the left hand surjection and then applying the clockwise sequence of maps. 
\end{lemma}
\begin{proof}
In order to see that the diagram commutes, given $[g]\in \pi_1(X^j)$ let $[\tilde{g}]\in\pi_1(X^k)$ be any choice of `lift' of $[g]$ via
$$
\xymatrix{
\pi_1(X^j, x_j) & \ar@{->>}[l]_{\pi_1E_0} \pi_1(H^{jk},\gamma^{jk}) \ar[r]^{\pi_1E_1}& \pi_1(X^k,x_k).
} 
$$
The assumption that $E_0^{-1}(x_j)\subset H^{jk}$ is connected ensures that the left hand map is surjective so that this is always possible. Then, possibly after replacing $g$ by another representative of its class in $\pi_1(X^j,x_j)$, there is a homotopy $\gamma^{jk}\tilde{g} \simeq g \gamma^{jk}$ (see Figure \ref{composition} below) so that for any $[h] \in \pi_1(X^k)$ we have
$$
\Gamma^{ik}[\tilde{g}h] = [\gamma^{ij}\gamma^{jk} \tilde{g}h]  =
[\gamma^{ij}g\gamma^{jk}h] = \Gamma^{ij}[g] \circ \Gamma^{jk}[h]
$$ 
as required.
\end{proof}

\begin{figure}[htbp]
\begin{center}
\centerline {
\includegraphics[width=3.5in]{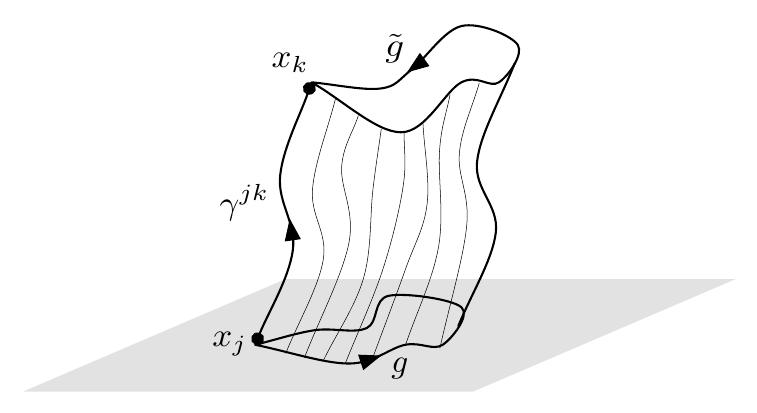}
}
\caption{Picture of the homotopy $\gamma^{jk}\tilde{g} \simeq g \gamma^{jk}$ from Lemma \ref{composing}.}
\label{composition}
\end{center}
\end{figure}

\section{Stratified covers}
\label{stratified covers}

Suppose $X$ is locally path-connected and locally simply-connected. Then the categories of covers of $X$, of spaces over $X$ with the unique lifting property for homotopies, of local systems on $X$ and of set-valued representations of the fundamental groupoid of $X$, \ie functors from the fundamental groupoid to sets, are equivalent. (Here we say that a map $p:Y\to X$ has the unique lifting property for homotopies if there is a unique solution to the lifting problem
$$
\xymatrix{
A\times \{0\} \ar[r] \ar@{_{(}->}[d] & Y \ar[d]^p \\
A \times [0,1] \ar[r] \ar@{-->}[ur] & X
}
$$
where $A$ is any CW-complex.) 

Now suppose $X$ is homotopically stratified. We will generalise these equivalences to the case when the strata of $X$ are locally path-connected and locally simply-connected. This condition is equivalent to asking that $X$ be `locally po-path-connected and locally po-simply-connected' in the following sense.
\begin{lemma}
\label{connectivity lemma}
Suppose $X$ is homotopically stratified. Then the strata of $X$ are locally path-connected and locally simply-connected if and only if each point $x$ of $X$ has a neighbourhood $U$ such that, up to homotopy through po-paths, there is a unique po-path from $x$ to any $x'\in U$, equivalently, if and only if $x$ is an initial object in $\Pi^{po}_1U$.
\end{lemma}
\begin{proof}
If each $x$ has a neighbourhood $U$ such that $x$ is initial in $\Pi^{po}_1U$ then it is clear that the strata must be locally path-connected and simply-connected. Suppose then that the strata have these properties. Each stratum $X^i$ is an almost stratum-preserving deformation retract of a neighbourhood $N$ of $X^i$ in $X$ \cite[Proposition 3.2]{quinn}. Let $U$ be the inverse image under this retraction of a connected and simply-connected neighbourhood of $x $ in $X^i$. Let $x' \in U$. We show that $\Pi^{po}_1U(x,x')$ has a unique element. Let $\rho$ be the path from a point of $X^i\cap U$ to $x'$ given by the retraction. Then there is a po-path from $x$ to $x'$: compose a path in $X^i\cap U$ from $x$ to $\rho(0)$ with $\rho$. To see that this po-path is unique up to homotopy, note that any po-path $\gamma$ from $x$ to $x'$ is homotopic through po-paths to the composition of a path in $X^i\cap U$ with $\rho$. Specifically $\gamma$ is homotopic to the result of applying the retraction to $\gamma$ composed with $\rho$. The result follows from the fact that $U\cap X^i$ is simply-connected.
   \end{proof}

In order to state our result we need to introduce appropriate generalisations of covers, local systems and so on. In each case there are two ways to relax the definition in a stratified context. For the remainder of this section, assume that $X$ is a homotopically stratified set and that the strata are locally path-connected and locally simply-connected. 

\begin{figure}[htbp] 
   \centering
  \includegraphics[width=2in]{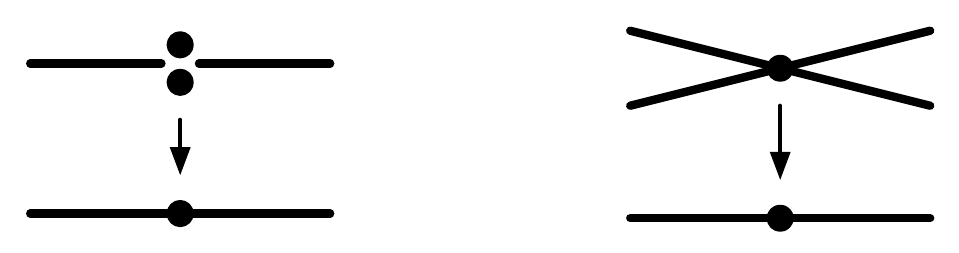} 
   \caption{The `line with two origins' is naturally a stratified \'etale cover of the real line stratified by $\{0\}$ and $\R - \{0\}$; two transversely intersecting lines a stratified branched cover. }
   \label{covers}
\end{figure}

The generalisations of covers are maps $p:Y\to X$ which restrict to covering maps over each stratum and are either
\begin{enumerate}
\item \'etale, \ie local homeomorphisms or
\item  locally-connected uniquely-complete spreads (see Appendix \ref{cosheaves}).
\end{enumerate}
The second class are branched covers in the topological sense, see Appendix \ref{cosheaves} for the definition. We refer to such maps as stratified \'etale covers and stratified branched covers respectively. Figure \ref{covers} illustrates motivating examples. Let $\et{X}$ and $\br{X}$ be the categories with respective objects the stratified \'etale and branched covers over $X$ and with maps the continuous maps over $X$.

More evidently, the two ways to relax the homotopy-theoretic notion of cover as a space with the unique homotopy lifting property are to consider maps $p:Y \to X$ which have the unique lifting property either
\begin{enumerate}
\item  for families of po-paths or
\item  for families of op-paths.
\end{enumerate}
Here by an op-path we mean the reverse $\overline{\gamma} : [0,1] \to X: t \mapsto \gamma(1-t)$ of a po-path. A \defn{family} of po-paths is a continuous map
$$
F:A \times [0,1] \to X 
$$
where $A$ is a CW-complex and for each $a\in A$ the restriction $F(a,-)$ is a po-path. Families of op-paths are defined analogously. Let $\ulpo{X}$ and $\ulop{X}$ respectively be the categories of spaces over $X$ with these properties. In both cases morphisms are continuous maps over $X$.

\begin{proposition}
\label{lifting}
A stratified \'etale cover has the unique lifting property for families of po-paths. A stratified branched cover has the unique lifting property for families of op-paths. 
\end{proposition}
\begin{proof}
We begin with the cases of a single po-path or op-path and then generalise to families. Note that if $\gamma : [0,1] \to X$ is a po-path and $\gamma(t)$ is in a stratum $X^i$ then there is some $\epsilon > 0$ such that $\gamma(s)\in X^i$ for $s\in (t-\epsilon, t]$. Obviously for an op-path there is a similar statement but with $\gamma(s)\in X^i$ for $s\in [t, t+\epsilon)$.

Suppose that $p:Y\to X$ is a stratified \'etale cover and that $\gamma : [0,1] \to X$ is a po-path. Let $\widetilde{\gamma}(0)$ be a lift of $\gamma(0)$ to $Y$ and let $L \subset [0,1]$ be the set of $s$ for which $\gamma|_{[0,s]}$ has a unique lift starting at $\widetilde{\gamma}(0)$. Clearly $0\in L$. Suppose $[0,t) \subset L$. Then by the first observation $\gamma(t)$ lies in the same stratum as $\gamma(s)$ for $s\in (t-\epsilon, t]$.  Since $p$ is a cover over each stratum the unique lifting property of covers shows that $t \in L$. Hence $L$ is closed. On the other hand, it $t\in L$ then $\widetilde{\gamma}(t)$ has a neighbourhood $U$ such that $p|_U: U \to p(U)$ is a homeomorphism onto an open neighbourhood of $\gamma(t)$. It follows that $L$ is open. Since $[0,1]$ is connected we have $L=[0,1]$ and we can  lift po-paths uniquely. 

Suppose that $p:Y \to X$ is a stratified branched cover and that  $\gamma : [0,1] \to X$ is an op-path. Let $\widetilde{\gamma}(0)$ be a lift of $\gamma(0)$ to $Y$ and let $L \subset [0,1]$ be the set of $s$ for which $\gamma|_{[0,s]}$ has a unique lift $\widetilde{\gamma}:[0,s]\to Y$ starting at $\widetilde{\gamma}(0)$. Clearly $0\in L$. Suppose $[0,t) \subset L$. Then the trace of $\widetilde{\gamma}$ determines an element of
$$
\lim_{U \ni \gamma(t)} \pi_0(p^{-1}U) \cong p^{-1}\gamma(t).
$$
This is the unique continuous extension $\widetilde{\gamma}(t)$. It follows that $L$ is closed. Now suppose $[0,t]\subset L$. Then by the observation at the beginning of the proof, $\gamma(s)$ is in the same stratum for $s\in [t,t+\epsilon)$ for some $\epsilon >0$. Since $p$ restricts to a covering of each stratum we can extend the lift uniquely and $L$ is open. Therefore $L=[0,1]$ and we can uniquely lift op-paths.

To deal with the family case it remains only to show that these unique lifts fit into continuous families. The continuity of the lift of a family of po-paths to a stratified \'etale cover follows easily from the local homeomorphism property of the cover. Thus we will focus on the case of lifting a family $F:A\times[0,1]\to X$ of op-paths, parameterised by a CW-complex $A$, to a stratified branched cover $p:Y \to X$. Lifting each op-path individually uniquely defines a lift $\widetilde{F}:A\times[0,1]\to Y$ which is certainly continuous on $\{a\}\times[0,1]$ for each $a\in A$. To show that it is continuous on $A\times [0,1]$ we use the fact that $f:Z\to Y$ is continuous at $z$
\begin{itemize}
\item[$\iff$] for each open neighbourhood $V\ni f(z)$ the inverse image $f^{-1}V$ contains an open neighbourhood of $z$,
\item[$\iff$] for each open neighbourhood $U\ni pf(z)$ there is an open neighbourhood of $z$ mapping to the component $V$ of $f(z)$ in $p^{-1}U$.
\end{itemize}
If $\widetilde{F}$ is continuous on $A\times[0,t)$ then we can use this together with  the continuity of $\widetilde{F}$ on each $\{a\}\times [0,1]$ and the local-connectivity of $A$ to show that $\widetilde{F}$ is continuous on $A\times[0,t]$. 

Now suppose $\widetilde{F}$ is continuous on $A\times[0,t]$. Fix $a\in A$. By the initial observation $F(a,s)$ is in the same stratum for $s\in [t,t+\delta]$ for some $\delta >0$. As $Y$ is a stratified branched cover the cosheaf of components of $Y$ (see appendix \ref{cosheaves}) is locally-constant on strata. This means that we can cover 
$$
\{ F(a,s)\ |\ s\in [t,t+\delta] \}
$$
by a finite sequence of open neighbourhoods $U_i \ni F(a,t_i)$ for $t=t_0< t_1< \cdots <t _n=t+\delta$ each having the property that the evident map
$$
p^{-1}x \to \pi_0\left(p^{-1}U_i\right)
$$
is an isomorphism for each $x\in U_i$. It follows that if $\widetilde{F}(a,t)$ and $\widetilde{F}(a',t)$ are in the same component of $p^{-1}U_0$ and $$F(a,s) , F(a',s) \in U_0\cup\cdots \cup U_n \qquad \forall s\in [t,t+\delta]$$ then $\widetilde{F}(a,t+\delta)$ and $\widetilde{F}(a',t+\delta)$ are in the same component of $p^{-1}U_n$. Thus continuity propagates from $(a,t)$ to $(a,t+\delta)$ and $\widetilde{F}$ is continuous on the whole of $A\times [0,1]$.
\end{proof}

It follows immediately from the definition that a space over $X$ with the unique lifting property for families of po-paths defines a functor $\Pi_1^{po}X \to \cat{Set}$. Similarly, a space over $X$ with the unique lifting property for families of op-paths defines a functor $\Pi_1^{op}X \to \cat{Set}$, where $\Pi_1^{op}X$ is the analogue of $\Pi_1^{po}X$ but with po-paths replaced by op-paths. We now explain how to construct generalisations of local systems from such functors.

A local system is a locally-constant sheaf, but it is also a locally-constant cosheaf (see Appendix \ref{cosheaves} for the definition). The two appropriate generalisations of a local system to the stratified context are constructible sheaves and constructible cosheaves. These are, respectively, sheaves and {\em spatial} cosheaves which are locally-constant when restricted to each stratum of $X$. See \ref{spatial} for the definition of spatial (pre)cosheaf. Denote the resulting full subcategories of sheaves and cosheaves by $\shc{X}$ and $\cshc{X}$.

Given a functor $F: \Pi_1^{po}X \to \cat{Set}$ we define a presheaf $\mathcal{F}$ on $X$  with $\mathcal{F}(U)$ being the set of functions  $f: U \to \prod_{x\in U} F(x)$ such that  $f(x)\in F(x)$ and $f(x') = F(\gamma)\left(f(x)\right)$ whenever $\gamma$ is a po-path in $U$ from $x$ to $x'$. It follows from Lemma \ref{connectivity lemma} that the stalk $\mathcal{F}_x =F(x)$. Furthermore, the restriction of $\mathcal{F}$ to a stratum is a locally-constant presheaf. Hence the sheafification is a constructible sheaf with stalk $F(x)$ at $x\in X$. 

Similarly given a functor $G: \Pi_1^{op}X \to \cat{Set}$ we let $\mathcal{G}$ be the precosheaf with cosections $\mathcal{G}(U) = \sum_{x\in U} G(x) \big/ \sim$ where $\sim$ is the equivalence relation generated by $a \sim a'$ if there is an op-path $\gamma$ in $U$ from $x$ to $x'$ with $a'=G(\gamma)(a)$. This is locally-constant on strata and, using Lemma \ref{connectivity lemma}, we see that the costalk $\mathcal{G}_x=G(x)$. By Lemma \ref{spatial example}, this precosheaf is actually a spatial cosheaf. Therefore it is a constructible cosheaf. 

Finally, we can construct an \'etale space over $X$ from any sheaf on $X$. If the sheaf is constructible then the \'etale space will be a cover over each stratum of $X$, \ie it will be a stratified \'etale cover. Less well-known is the fact (see Appendix \ref{cosheaves}) that we can construct a locally-connected, uniquely-complete spread from a spatial cosheaf. In particular, when the cosheaf is constructible the corresponding spread is a stratified branched cover.

\begin{theorem}
\label{representations}
There are  equivalences of categories 
$$
\et{X} \simeq \ulpo{X} \simeq \funct{\Pi_1^{po}(X)}{\cat{Set}} \simeq \shc{X}
$$
and
$\br{X} \simeq\ulop{X} \simeq \funct{\Pi_1^{op}(X)}{\cat{Set}} \simeq \cshc{X}$. 
\end{theorem}
\begin{proof}
We have shown that there are maps of the objects
$$
\xymatrix{
\et{X} \ar[r] & \ulpo{X} \ar[d] \\
\shc{X} \ar[u] & \funct{\Pi_1^{po}(X)}{\cat{Set}} \ar[l]
}
\quad \textrm{and} \quad
\xymatrix{
\br{X} \ar[r] & \ulop{X} \ar[d] \\
\cshc{X} \ar[u] & \funct{\Pi_1^{op}(X)}{\cat{Set}}. \ar[l]
}
$$
It is tedious but easy to check that each construction is functorial and that starting at any category and cycling around the diagram is an auto-equivalence. We leave this as an exercise!
\end{proof}

\begin{remark}
\label{rep rem}
If a functor $F:\cat{C}\to \cat{D}$ induces an equivalence $F^*:\funct{\cat{D}}{\cat{Set}} \to \funct{\cat{C}}{\cat{Set}}$ then $F$ is fully faithful. Therefore $F$ is an equivalence if, in addition, it is essentially surjective. To see that $F$ is fully faithful, note that $F^*$ has a left adjoint $L$ (which is also fully faithful). It follows from the Yoneda lemma that $L$ takes the representable functor  $\hom_C(c,-)$ to the representable $\hom_D(Fc,-)$, and that it takes composition with $f:c' \to c$ to composition with $Ff: Fc'\to Fc$, \ie that the diagram
$$
\xymatrix{
\cat{C}^{op} \ar[r]^{F^{op}} \ar[d] & \cat{D}^{op} \ar[d]\\
\funct{\cat{C}}{\cat{Set}} \ar[r]_{L} & \funct{\cat{D}}{\cat{Set}},
}
$$
in which the vertical functors are the Yoneda embeddings, commutes. The claim follows.
\end{remark}

\begin{corollary}
\label{ho to po}
If $X$ is homotopically stratified and each stratum is locally path-connected and locally simply-connected then the functor $\Pi^{ho}_1X \to \Pi^{po}_1X$ arising from the inclusion of the homotopy link in the space of po-paths is an equivalence.
\end{corollary}
\begin{proof}
Theroem \ref{representations} and the preceding discussion can be repeated with $\Pi^{ho}_1X$ in place of $\Pi^{po}_1X$. We conclude that $\Pi^{ho}_1X \to \Pi^{po}_1X$ induces an equivalence 
$$
\funct{\Pi^{po}_1X}{\cat{Set}} \to \funct{\Pi^{ho}_1X}{\cat{Set}},
$$
and so is an equivalence by Remark \ref{rep rem}.
\end{proof}

\begin{corollary}
\label{fund groupoid}
Suppose $X$ is a  homotopically stratified space with locally-connected and locally simply-connected strata.  Let $\Pi^{po}_1X_{loc}$ be the category obtained by localising the fundamental category at the set of all morphisms --- thus the objects of $\Pi^{po}_1X_{loc}$ are the points of $X$ and the morphisms are equivalence classes of words in the morphisms of $\Pi^{po}_1X$ and formal inverses thereof under obvious relations arising from composition and cancellation. The functor $\Pi^{po}_1X_{loc} \to \Pi_1X$ given by composing the terms in these words is an equivalence, \ie $\Pi_1X$ is the `groupoidification' of $\Pi^{po}_1X$.
\end{corollary}
\begin{proof}
The obvious functors $\Pi_1^{po}X \to \Pi_1^{po}X_{loc} \to \Pi_1X$ induce functors
$$
\funct{\Pi_1X}{\cat{Set}} \to \funct{\Pi^{po}_1X_{loc}}{\cat{Set}} \to \funct{\Pi^{po}_1X}{\cat{Set}}.
$$
It follows from the connectivity assumptions on the strata that $X$ is locally-connected and locally simply-connected. Hence  $\funct{\Pi_1X}{\cat{Set}}$ is equivalent to the category of covers of $X$. On the other hand $\funct{\Pi^{po}_1X_{loc}}{\cat{Set}} \to \funct{\Pi^{po}_1X}{\cat{Set}}$ is the inclusion of the full subcategory of functors which take all maps to isomorphisms. A stratified \'etale cover is a genuine cover precisely when the monodromy induced by any po-path is an isomorphism. Thus  $$\funct{\Pi^{po}_1X_{loc}}{\cat{Set}} \to \funct{\Pi^{po}_1X}{\cat{Set}}$$ corresponds to the inclusion of the category of covers in the category of stratified \'etale covers. Hence $\funct{\Pi_1X}{\cat{Set}} \to \funct{\Pi^{po}_1X_{loc}}{\cat{Set}}$ is an equivalence. It follows from Remark \ref{rep rem} that $\Pi_1^{po}X_{loc} \to \Pi_1X$ is an equivalence.
\end{proof}

Thus the fundamental category of a homotopically stratified set (with locally-connected and simply-connected strata) contains at least as much information as the fundamental groupoid. In some cases this fact can be used to compute the fundamental group by computing a skeleton of $\Pi^{po}_1 X$, localising this to obtain a groupoid --- which is equivalent to the fundamental groupoid by the above theorem --- and then reading off the fundamental group as the automorphisms of an object. The following example is offered as a `proof of concept', not because it is an elegant way to compute $\pi_1\R\PP^n$!
\begin{example}
\label{rp2fg}
Real projective space $\R\PP^n$ has a filtration $\R\PP^0 \subset \R\PP^1 \subset \cdots \subset \R\PP^n$. Let $f:S^n \to \R\PP^n$ be the standard $2$ to $1$ covering, and $g:\R^{n+1} -\{0\} \to \R\PP^n$ the standard quotient map. Let 
$$
y_i = (0,\ldots, 0, 1,0, \ldots,0)
$$
where the non-zero entry is in the $i$th place and choose the basepoint $x_i=f(y_i)$ for the stratum $X^i := \R\PP^i - \R\PP^{i-1}$. The full subcategory on the objects $x_i$ is a skeleton of $\Pi^{po}_1 \R\PP^n$ because every point is connected to precisely one of the $x_i$ by a reversible po-path.  In order to characterise po-paths it is helpful to define the function 
$$
L : \R\PP^n \to \{0,\ldots,n\} : [p_0 : \ldots : p_n] \mapsto \max \{j \ | \ p_j\neq 0\}.
$$
Then $p\in X^i \iff L(p)=i$, and a path $\gamma: [0,1] \to \R\PP^n$ is a po-path if and only if the composite $L\circ \gamma : [0,1] \to \{0,\ldots,n\}$ is increasing.

Let $\gamma : [0,1]\to \R\PP^n$ be a po-path from $x_i$ to $x_j$ where $i< j$. (The case $i=j$ is uninteresting because the strata are simply-connected.) Let $\tilde\gamma : [0,1] \to S^n$ be the unique lift of $\gamma$ along the covering map $f$ starting at $y_i$. The end point of $\tilde \gamma$ is then at  $\epsilon_\gamma y_j$ where $\epsilon_\gamma \in \{\pm 1\}$. Furthermore, by homotopy lifting, the end point is the same for any po-path homotopic to $\gamma$ through po-paths so there is a well-defined map
$$
\epsilon : \Pi^{po}_1\R\PP^n(x_i,x_j) \to \{\pm 1\}.
$$
In fact this is a bijection. To see this consider the map
$$
\tilde\eta:[0,1]^2 \to \R^{n+1} -\{0\} : (s,t) \mapsto (1-s) \tilde \gamma(t) + s(0,\ldots,1-t,\ldots,\epsilon_\gamma t,\ldots,0)
$$
where the only non-zero entries of the right hand term are in the $i$th and $j$th places. The composite $\eta = g\circ \tilde \eta$ is a homotopy relative to end points from $\gamma$ to an element of  $\holink{X^i\cup X^j}{X^i}$. To see that it is a homotopy through po-paths note that, because $\gamma$ is a po-path, we have $\epsilon_\gamma \tilde \gamma_j(t) \geq 0$ for $t\in [0,1]$. It follows that $L(\eta(s,t))$ is an increasing function of $t$ for each $s\in [0,1]$ as required.

Hence $\Pi^{po}_1\R\PP^n(x_i,x_j) = \{\alpha,\beta\}$ is a two element set. (Properly we should write $\alpha_{ij}$ and $\beta_{ij}$ but we omit the subscripts for ease of reading.) It is easy to check that $\epsilon_{\gamma \cdot \gamma'} = \epsilon_\gamma \epsilon_{\gamma'}$ so that  composing paths we have
$\alpha^2 = \beta^2$ and $\alpha \beta = \beta \alpha$. Localising introduces inverses $\alpha^{-1}$ and $\beta^{-1}$ satisfying $$\alpha\beta^{-1} = \alpha^{-1}\beta = \beta\alpha^{-1}=\beta^{-1}\alpha.$$  The automorphisms of $x_0$ in the resulting groupoid  are 
$$
\langle \alpha\beta^{-1}\  |\  (\alpha\beta^{-1})^2=1 \rangle \cong \Z / 2\Z \cong \pi_1(\R\PP^n,x_0)
$$
as expected from Corollary \ref{fund groupoid}.
\end{example}

\section{Configuration spaces of points in the plane}
\label{config}

In this section we compute (a skeleton of) $\Pi^{po}_1X$ when $X$ is the configuration space of $n$ indistinguishable but not necessarily distinct points in $\C$, \ie $X$ is the symmetric product $SP^n\C:= \C^n / S_n$ where the symmetric group $S_n$ acts on $\C^n$ by permuting coordinates. A point $x$ of $SP^n\C$ determines a configuration of $n$ indistinguishable points in $\C$ given by the set of coordinates of any pre-image of $x$ in $\C^n$.

The symmetric product has a natural stratification by orbit types indexed by partitions of $n$. The strata are the subsets of points corresponding to configurations where the $n$ points coalesce according to the indexing partition. More precisely, let $C$ be a concrete partition of $\{1,\ldots,n\}$, \ie $C$ is a set $\{C_1,\ldots,C_k\}$ of disjoint subsets of $\{1,\ldots,n\}$ whose union is the entire set. Define
$$
Y^C = \{ (z_1,\ldots,z_n\} \in \C^n \ | \ z_a=z_b \iff a,b \in C_i \ \textrm{for some}\ 1\leq i\leq k\}.
$$
Then $Y^C$ is a closed complex submanifold of an open subset of $\C^n$ and has codimension
$$
\sum_{i=1}^k \left( |C_i| -1 \right).
$$
It is obtained by deleting complex linear subspaces from a complex linear subspace and is therefore connected. For a concrete partition $C$ let $\mathcal{P}(C)$ be the corresponding (abstract) partition of $n$ with cardinalities $|C_1|, \ldots, |C_k|$. For a partition $P$ set
$$
Y^P = \bigcup_{\mathcal{P}(C)=P} Y^C.
$$
Note that there is an element $\pi \in S_n$ inducing a complex-analytic isomorphism $Y^C \cong Y^{C'}$ exactly when $\mathcal{P}(C)=\mathcal{P}(C')$. Hence the $Y^P$ are invariant under the action of $S_n$ on $\C^n$. The quotients $X^P=Y^P/S_n$ form a stratification of $X$ by connected complex analytic strata.

The poset of strata is the set of partitions with the relation $P \leq Q$ if $Q$ is a refinement of $P$, that is if $Q$ is obtained by further partitioning the parts of $P$. We denote the number of parts in $P$ by $|P|$ and write the partition with $k$ parts of cardinalities $p_1,\ldots,p_k$ as $(p_1|\cdots|p_k)$. In this notation, the top element, corresponding to the open stratum, is $(1 | 1 | \cdots |1)$ and the bottom, corresponding to the stratum where all $n$ points coalesce, is $(n)$.

\begin{figure}[htbp] 
   \centering
   \includegraphics[width=2.5in]{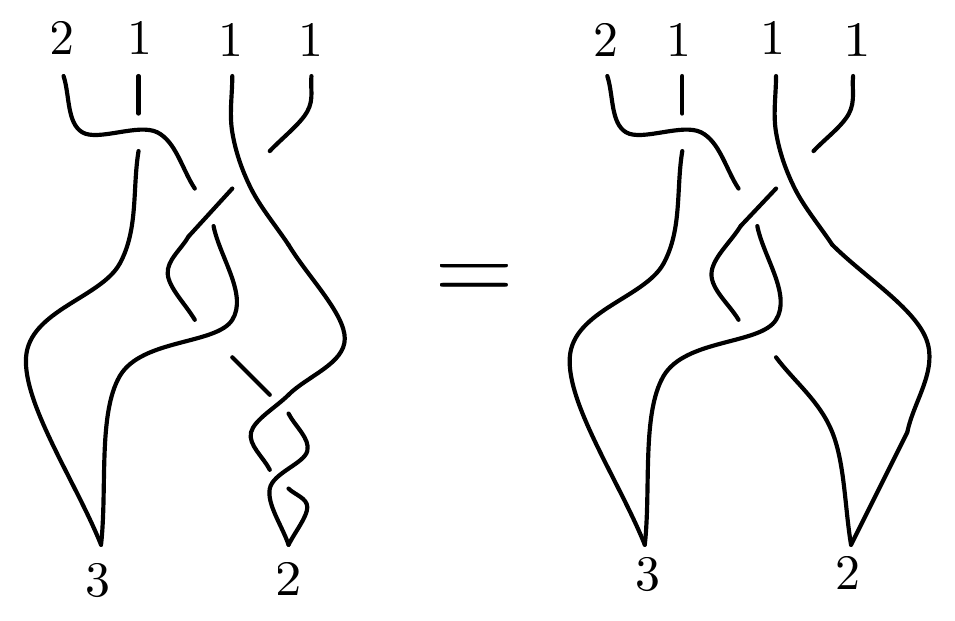} 
   \caption{Two representatives of a morphism in the fundamental category from $x_P$ to $x_Q$ where $P=(3|2)$ and $Q=(2|1|1|1)$.}
   \label{glued braid}
\end{figure}

Fix a basepoint $x_P$ in each stratum.  Elements of $\Pi^{po}_1X(x_P,x_Q)$ are represented by po-paths from $x_P$ to $x_Q$ in the homotopy link. Thinking of points of the symmetric product as configurations in $\C$, the graph of such a po-path corresponds to a set of strings in $\C \times [0,1]$ joining a configuration of type $P$ to one of type $Q$. The condition that it is a po-path in the homotopy link means that the set of strings is a braid on $|Q|$ strings where the starts of the strings are glued together according to the partition $P$. Furthermore, each string is labelled by a natural number which records the cardinality of the corresponding part of $Q$. The sum of these for the set of strings emanating from a point, which corresponds to a part of $P$, is the cardinality of that part. Two representations give the same morphism if they are isotopic relative to the end points. Thus there are the usual braid relations on $|Q|$ strings, but also relations coming from  `internal' braiding within the parts of $P$ which becomes `external' braiding of the $|Q|$ strings.  Figure \ref{glued braid} illustrates this in a simple example. 

Let's make this precise. In order that we may work with subgroups of $S_n$ and of the braid group $B_n$ we choose a po-path $\gamma^{P,Q}$ from $x_P$ to $x_Q$ whenever $P<Q$. We do this in such a way that $\gamma^{P,R}$ is the composite path $\gamma^{P,Q} \cdot \gamma^{Q,R}$ whenever $P < Q< R$. Label the points in the configuration corresponding to the basepoint in the open stratum by $1,\ldots,n$. The chosen po-paths then identify  a concrete partition $C$ of $\{1,\ldots,n\}$ with $\mathcal{P}(C)=P$ corresponding to each basepoint $x_P$: namely $i$ and $j$ are in the same part at $x_P$ if the points $i$ and $j$ coalesce along the reverse of the po-path from $x_P$ to $x_{(1|\cdots|1)}$. Henceforth we will abuse notation by denoting this choice of concrete partition by the same letter as the corresponding abstract partition.

To each concrete partition $P$ of $\{1,\ldots,n\}$ into $k$ subsets of cardinality $p_1,\ldots,p_k$ we associate three subgroups of $S_n$ and three subgroups of $B_n$. Let
 $$
 S_P = \{ \sigma\in S_n : i \sim_P j \iff \sigma(i) \sim_P \sigma(j)\},
 $$ 
 where $\sim_P$ is the equivalence relation on $\{1,\ldots,n\}$ corresponding to $P$. This group of symmetries of $P$ is a product:
$$
S_P \cong IS_P \times ES_P
$$
where the internal symmetries $IS_P\cong S_{p_1} \times \cdots \times S_{p_k}$ permute the elements within the parts and the external symmetries $ES_P \leq S_k$ permute the parts amongst themselves, preserving the cardinality. Similarly we define $B_P$ to be the product $IB_P \times EB_P$ of the internal braids $IB_P \cong B_{p_1} \times \cdots \times B_{p_k}$, which braid the strings within the parts, and the external braids $EB_P$. The latter are defined by the pullback square
$$
\xymatrix{
EB_P \ar[r] \ar[d] & B_k \ar[d] \\
ES_P \ar[r] & S_k.
}
$$
The partition $P$ defines embeddings of $IB_P$ and $EB_P$ into $B_n$ as commuting subgroups. This determines an embedding of $B_P$ into $B_n$.

It is well-known --- and visually obvious --- that the fundamental group of the top stratum $X_{(1|\cdots|1)}$  is the braid group $B_n$ on $n$ strings. More generally, for a partition $P$ into subsets of cardinality $p_1,\ldots,p_k$, the fundamental group of the stratum $X^P$ is isomorphic to $EB_P$. We do not get the full braid group on $k$ strings because points of the stratum correspond to sets of $k$ points {\em labelled by the cardinalities $p_1,\ldots,p_k$} and this labelling must be preserved by the braiding. There are also no `internal' braids.

Pick partitions $P \leq Q$. Let $E_0:\holink{X^P\cup X^Q}{X^P}\to X^P$ be the evaluation at $0$ map. From Figure \ref{braided LES} there is an exact sequence
\begin{equation}
\label{config es}
\pi_1\left(E_0^{-1}x_P, \gamma^{P,Q}\right) \to \pi_1\left(X^Q,x_Q\right) \to \Pi^{po}_1(x_P,x_Q) \to \pi_0\left(E_0^{-1}x_P\right) \to 1.
\end{equation}
The fibre $E_0^{-1}x_P$ is the space of `geometric braids' embedded in $\C\times[0,1]$ starting at the configuration $x_P$ and immediately splitting into $|Q|$ strings. These strings are labelled by the cardinalities of the parts of $Q$ in such a way that the sum of the labels on the set of strings emanating from the same point of $x_P$ is the cardinality of the corresponding part of $P$.

Since the ends are free to move the set of components $\pi_0\left(E_0^{-1}x_P\right)$ is the set of ways in which $P$ can be refined to $Q$. We can describe this in terms of concrete partitions of $\{1,\ldots,n\}$ as follows. Define a sub\emph{set}
$$
S_{P,Q} = \{ \sigma \in S_n : \sigma(i) \sim_Q \sigma(j) \Rightarrow i \sim_P j \}.
$$
The subgroup $IS_P$ acts on the left and $S_Q$ acts on the right. The set of ways of refining $P$ to $Q$ is the double orbit space:
$$
\pi_0\left(E_0^{-1}x_P\right) \cong IS_P \backslash S_{P,Q} / S_Q.
$$ 

The image of $\pi_1E_1:\pi_1\left(E_0^{-1}x_P, \gamma^{P,Q}\right) \to \pi_1\left(X^Q,x_Q\right) \cong EB_Q$ consists of braids on $|Q|$ strings (respecting labels) which can be continuously deformed by gathering together strings in the same part of $P$ so that they become the trivial braid on $|P|$ strings. Thus the image is precisely the subgroup $IB_P \cap EB_Q$ (where the intersection is taken by embedding both in $B_n$).

Therefore (\ref{config es}) gives an exact sequence
\begin{equation}
\label{config es 2}
1\to IB_P \cap EB_Q \backslash EB_Q \to \Pi^{po}_1(x_P,x_Q) \to IS_P \backslash S_{P,Q} / S_Q \to 1.
\end{equation}
Here is an explicit way to construct the middle term. Let $EB_{P,Q}$ be the sub\emph{set} of $B_{|Q|}$ defined by the pullback square
$$
\xymatrix{
EB_{P,Q} \ar[r] \ar[d] & B_{|Q|} \ar[d] \\
S_{P,Q} \ar[r] & S_n
}
$$
where the map $B_{|Q|} \to S_n$ is obtained by permuting the parts of the partition at $x_Q$ whilst leaving the ordering within parts fixed. Note that $S_{Q,Q}=S_Q$, and also $EB_{Q,Q} = EB_Q$ because $EB_Q$ fits into the pullback square
$$
\xymatrix{
EB_Q \ar[r] \ar[d] & B_{|Q|} \ar[d] \\
S_Q \ar[r] & S_n.
}
$$
It follows from this and the fact that $S_Q \subset S_{P,Q}$ that there is an inclusion $EB_Q \to EB_{P,Q}$. Furthermore, the group $S_Q$ acts on the right on the set $S_{P,Q}$ and it follows that $EB_Q$ acts on $EB_{P,Q}$ with orbit space $S_{P,Q} / S_Q$. I.e., appropriately interpreted, there is an exact sequence
$$
1\to EB_Q \to EB_{P,Q} \to S_{P,Q}/S_Q \to 1.
$$
The partition $Q$ determines an embedding of $EB_{P,Q}$ into $B_n$, and henceforth we identify $EB_{P,Q}$ with its image in $B_n$. Denote the image of $EB_{P,Q}$ in the orbit space $IB_P \backslash B_n$ by $IB_P \backslash EB_{P,Q}$. Comparing with (\ref{config es 2}) it follows that
\begin{equation}
\label{hom set}
 \Pi^{po}_1(x_P,x_Q) \cong IB_P \backslash EB_{P,Q}.
\end{equation}

The composition in $\Pi^{po}_1X$ is simple to describe in terms of the composition in $B_n$. Suppose $P \leq Q \leq R$. Composition within $S_n$ determines a map $S_{P,Q} \times S_{Q,R} \to S_{P,R}$ and it follows that composition in $B_n$ gives a map
$$
EB_{P,Q} \times EB_{Q,R} \to EB_{P,R}.
$$
Noting that the internal braids $IB_Q$ commute with $EB_{P,Q}$ and are a subgroup of the internal braids $IB_P$ we see that the above map descends to 
$$
IB_P\backslash EB_{P,Q} \times IB_Q \backslash EB_{Q,R} \to IB_R \backslash EB_{P,R} : \left([\alpha],[\beta]\right) \mapsto [\alpha\beta].
$$
This is the composition in $\Pi^{po}_1X$.

The quotient map $\C^n \to \C^n/S_n = X$ is a stratified branched cover. It corresponds to the functor 
$$
x_P \mapsto IS_P \backslash S_n \quad \textrm{and} \quad\Pi^{po}_1(x_P,x_Q) \cong  IB_P \backslash EB_{P,Q} \mapsto IS_P \backslash S_{P,Q}
$$
where the map $IS_Q \backslash S_n \to IS_P \backslash S_n$ corresponding to an element of $IS_P \backslash S_{P,Q}$ is induced by composition $S_{P,Q} \times S_n \to S_n$ in $S_n$.

\appendix

\section{Pre-orders and spaces}
\label{pre-orders}
\label{alexandrov}

A pre-order $P$ is a set equipped with a reflexive and transitive relation $\leq$. An increasing map $f:P\to Q$ is one such that $f(p) \leq f(p')$ whenever $p\leq p'$. Let $\cat{Preorder}$ be the category of pre-orders and increasing maps. 

If $\leq$ is also antisymmetric, \ie $p \leq q$ and $q\leq p$ implies $p=q$, then $P$ is a poset. Each pre-order $P$ has an associated poset given by quotienting by the  equivalence relation $$p\sim p' \iff p\leq p'\ \textrm{and}\ p'\leq p.$$ The quotient map $P \to P/\!\sim$ is increasing. This construction is right adjoint to the natural inclusion $\cat{Poset}\hookrightarrow \cat{Preorder}$. 

Pre-orders arise naturally from topology: if $X$ is a topological space there is a pre-order on the points of $X$ given by 
$$
x \leq y \iff \left( x\in U \Rightarrow y\in U\right) \iff x\in \overline{y}
$$
where $\overline{y}$ is the closure of $y$, and $U$ is an open set. This is called the \defn{specialisation} pre-order because lower points are more `special' and higher ones more `generic' (in the sense familiar to algebraic geometers). It defines a functor $S: \cat{Top}\to \cat{Preorder}$. 

Conversely, given a pre-order $P$ we can topologise it (not necessarily uniquely) so that the resulting specialisation pre-order is $P$. The coarsest topology with this property has closed sets generated (under finite union and arbitrary intersection) by the $D_p=\{q\ |\ q\leq p\}$ for $p\in P$. The finest topology with this property is the \defn{Alexandrov topology}\label{alex top}, in which the sets $U_p = \{q \ |\ p \leq q\}$ form a basis. Write $D(P)$ for $P$ with the coarsest topology and $U(P)$ for $P$ with the Alexandrov topology. Here $D$ stands for `downward-closed' and $U$ for `upward-open'. The identity maps
$$
US(X) \to X \to DS(X)
$$
are continuous, so $U$ is a candidate left adjoint and $D$ a candidate right adjoint for $S$. It turns out that $D$ is not a right adjoint, indeed it is not even a functor; an increasing map $P\to Q$ need not induce a continuous map $D(P) \to D(Q)$. For example if we order $[0,1]$ in the usual way then 
$$
t \mapsto \left\{
\begin{array}{ll}
 t/2 &  t \in [0,1)   \\
 1 & t=1 
\end{array}
\right.
$$
is increasing but not continuous as a map $D[0,1] \to D[0,1]$. However, a map of pre-orders is increasing if, and only if, it is continuous in the Alexandrov topologies. In particular 
$$
U: \cat{Preorder} \to \cat{Top}
$$
is a functor and is left adjoint to specialisation. The identity maps on the underlying sets give the unit and counit of the adjunction
$$
P \to SU(P) \quad \textrm{and} \quad US(X)\to X.
$$
It follows from the equivalences $p\leq q \iff q \in U_p \iff U_q \subset U_p$ that $P\cong SU(P)$. However, the topology of $US(X)$ can be much finer than that of $X$. For instance if $X$ is a metric space then $US(X)$ has the discrete topology. In fact, $US(X)\cong X$ if and only if $X\cong U(P)$ for some pre-order $P$, in which case we say $X$ is an \defn{Alexandrov space}. 

Alexandrov spaces can be alternatively characterised as those spaces for which each point $x$ has a unique minimal open neighbourhood $U_x$ (by minimal we mean that for any open $U$ we have $x\in U \iff U_x\subset U$). A simple consequence is that any space with a finite topology is Alexandrov. Finally note that, if $P$ is a pre-order whose associated poset is finite then $U(P)\cong D(P)$ and there is a unique topological space whose specialisation pre-order is $P$.

\subsection{Pre-ordered and filtered spaces}
\label{filtered}

A pre-ordered space, or \defn{po-space} for short, is a pair $(X,\leq)$ consisting of a topological space $X$ and a pre-order $\leq $ on the points of $X$. A map of po-spaces, or \defn{po-map}, $f:X\to Y$ is a continuous and increasing map. Let $\cat{PoSpace}$ be the resulting category.

We can think of po-spaces in two other ways. The first is as a space equipped with two topologies, a `spatial' topology and an Alexandrov topology which defines the pre-order. Po-maps correspond to maps which are continuous with respect to the spatial and Alexandrov topologies. The second is as a space over a poset, \ie as a space $X$ together with a surjective map $\sigma_X: X\to P_X$ where $P_X$ is the poset associated to the pre-order on the points of $X$ and $\sigma_X$ the quotient map. We call $P_X$ the \defn{poset of strata} of $X$, for reasons which will become apparent in a moment. In this picture po-maps are commutative squares
$$
\xymatrix{
X \ar[d]_{\sigma_X} \ar[r]^f &Y \ar[d]^{\sigma_Y}\\
P_X \ar[r]_g & P_Y
}
$$
in which $f$ is continuous and $g$ increasing.

The definition of po-space assumes no compatibility between the topology and the pre-order. We now consider two compatibility conditions
\begin{description}
\item[C1]: the down-sets $D_x = \{ x'\ |\ x'\leq x\}$ are closed for all $x\in X$;
\item[C2]: the up-sets $U_x = \{ x' \ | \ x'\geq x\}$ are open for all $x\in X$.
\end{description}
The first of these is perhaps the most natural, we expect $\leq$ to be a closed condition. The second condition implies the first.

Po-spaces satisfying C1 are better known as filtered spaces: a topological space $X$ is \defn{filtered} if there are non-empty closed subspaces $X_i$ indexed by a poset $P_X$ such that 
\begin{enumerate}
\item $X_i \subset X_j \iff i\leq j$;
\item the subsets  $X^i = X_i - \bigcup_{j<i} X_j$ partition $X$.
\end{enumerate}
The $X^i$ are known as the \defn{strata} of the filtration and $P_X$  is referred to as the \defn{poset of strata}. There is an induced pre-order on the points of a filtered space $X$ coming from the map $\sigma_X: X \to P_X$ taking points to the stratum in which they lie. I.e.\ we define $x\leq x' \iff \sigma_X(x)\leq \sigma_X(x')$. By definition, for $x\in X_i$
$$
D_x = \{ x' \ | \ x'\leq x\} =  X_i
$$ 
is closed, so that a filtered space is a po-space satisfying C1. Conversely, if $(X,\leq)$ is a po-space satisfying C1 then we can filter it by the the non-empty closed subsets  $X_i = \{x' \ |\ x'\leq x\}$. 

The second compatibility condition C2 is equivalent, from the filtered perspective, to asking that upward unions of strata
$$
\bigcup_{j \geq i} X^i = X- \bigcup_{k \not\geq i} X_k
$$
are open for each $i$. Equivalently, the surjection $\sigma_X : X \to P_X$ is continuous with respect to the Alexandrov topology on $P_X$. We will say $X$ is \defn{well-filtered} when this holds --- it is automatic for a space with a finite filtration. The strata of a well-filtered space are locally-closed, \ie each is a closed subset of an open subset of the space.

\begin{example}
\label{interval1}
We define a well-filtered space $I_n$ whose underlying space is the interval $[0,1]$ with the standard topology and whose filtration is
$$
\{0\} \subset \left[0,\frac{1}{n}\right] \subset \cdots \subset \left[0,\frac{n-1}{n}\right]\subset [0,1].
$$
The strata are the subsets $\{0\}, (0,\frac{1}{n}], \ldots, (\frac{n-1}{n},1]$ and the poset of strata is $0\leq 1\leq  \cdots \leq n$. The Alexandrov topology on this has open sets 
$$
\{0,1,\ldots,n\}, \{1,2,\ldots,n\},\ldots,\{n\}\ \textrm{and}\ \emptyset.
$$ The continuous surjection corresponding to the filtration is given by $t \mapsto \lceil nt \rceil$. 

The ordered interval $I$ is the space $[0,1]$ with the standard topology and order. It is filtered because $D_t=\{s | \ s\leq t\} = [0,t]$ is closed, but not well-filtered because $U_t = [t,1]$ is not open for $t\neq 0$. 
\end{example}

A continuous map $f:X\to Y$ of filtered spaces is \defn{filtered} if there is an increasing map of the indexing posets $g: P_X \to P_Y$ such that $f(X_i) \subset Y_{g(i)}$.  This is a rather weak condition and we will instead consider the stronger notion of a \defn{stratified} map, \ie a map for which 
$$
f(X^i) \subset Y^{g(i)}.
$$ 
Every stratified map is filtered. The map $g$ can be recovered from $f$ (but the requirement that $g$ be increasing is a restriction on $f$). Note that $f:X\to Y$ is stratified if and only if
$$
\xymatrix{
X \ar[d]_{\sigma_X} \ar[r]^f &Y \ar[d]^{\sigma_Y}\\
P_X \ar[r]_g & P_Y
}
$$
commutes so that stratified maps are po-maps and vice versa. Hence the two compatibility conditions C1 and C2 between the topology and the pre-order cut out two full subcategories
$$
 \cat{PoSpace} \supset \cat{Filt} \supset \cat{WellFilt}
$$
consisting respectively of the filtered and the well-filtered spaces, and the stratified maps between them. Most interesting examples of po-spaces arising `in nature' seem to be filtered or even well-filtered, certainly we will focus on these. However, for the purposes of theory it is convenient to work with po-spaces.

\section{Cosheaves and complete spreads}
\label{cosheaves}

The notion of a sheaf and the correspondence of sheaves and \'etale spaces are well-known. In contrast cosheaves are rarely discussed and there are few references. A further complication is that, unlike the case of sheaves, a cosheaf of abelian groups is not simply a cosheaf of sets whose cosections have compatible abelian group structures. Thus the theories of cosheaves of sets and of abelian groups are different. (The reason is that, whilst the underlying set of a product of abelian groups is just the product of the underlying sets of the abelian groups, the same is not true for coproducts --- the coproduct of sets is the disjoint union whereas the coproduct of abelian groups is the direct sum.) 

This appendix provides the necessary background for the use of cosheaves in this paper, and in particular explains the correspondence between certain cosheaves, namely the spatial cosheaves, on a space $X$ and locally-connected, uniquely complete spreads over $X$. It is a simplified exposition of some of the results from \cite[\S 5,6]{funk}.

A \defn{precosheaf} of sets on a topological space $X$ is a functor $\cosh{F}: \open{X} \to \cat{Set}$ from the category of open subsets of $X$ and inclusions to the category of sets. Elements of $\cosh{F}(U)$ are called \defn{cosections} over $U$, and the maps $\cosh{F}(U)\to \cosh{F}(V)$ for $U \subset V$ are called \defn{extensions}. A \defn{cosheaf} of sets on $X$ is a precosheaf which preserves colimits, \ie for any collection $\{U_i\}$ of open sets $\cosh{F}\left(\bigcup_i U_i\right)$ is the colimit of 
$$
\sum_{i,j}\cosh{F}\left( U_i \cap U_j\right) \to \sum_i \cosh{F}\left( U_i\right).
$$
The displayed map is induced from the inclusions of $U_i \cap U_j$ into $U_i$ and $U_j$ in the obvious way. To be concrete, this means that
$$
\cosh{F}\left(\bigcup_i U_i\right) = \sum_i \cosh{F}\left( U_i\right) \big/ \sim
$$
is the quotient of the disjoint union by the equivalence relation generated by $\alpha_i \sim \alpha_j$ if there is $\beta \in \cosh{F}\left( U_i \cap U_j\right)$ which extends to both $\alpha_i \in \cosh{F}\left( U_i\right)$ and $\alpha_j \in \cosh{F}\left( U_j\right)$. Maps of precosheaves and cosheaves are natural transformations. We denote the categories of precosheaves and cosheaves on $X$ by $\cat{Precosh\!\downarrow\!X}$ and $\cat{Cosh\!\downarrow\!X}$ respectively.
\begin{example}
\label{component cosheaf}
Given a continuous map $p_Y:Y \to X$ the assignment
$$
U \mapsto \pi_0\left( p_Y^{-1}U \right)
$$
is a precosheaf $CY$. (Here $\pi_0$ denotes components, not path components.) This defines a functor $C: \cat{Top \!\downarrow\! X} \to \cat{Precosh\!\downarrow\!X}$. 

If $Y$ is locally-connected then $CY$ is in fact a cosheaf. To see this, recall that a space is locally-connected if and only if the connected components of open subsets are open. It follows that
$$
\pi_0 : \cat{LCTop} \to \cat{Set},
$$
where $\cat{LCTop}$ is the full subcategory of locally-connected spaces, is left adjoint to the discrete space functor. Since left adjoints preserve colimits $CY$ is a cosheaf when $Y$ is locally-connected; we call it the \defn{cosheaf of components} of $Y$. 
\end{example}

This example shows that we can naturally turn spaces over a given base into precosheaves on that space. Conversely, every precosheaf $\cosh{F}$ has an associated \defn{display space} $D\cosh{F} \in \cat{Top \!\downarrow\! X}$. As a set the display space is the disjoint union
$$
\sum_{x\in X} \cosh{F}_x
$$
where $\cosh{F}_x = \lim_{U\ni x} \cosh{F}(U)$ is the \defn{costalk} of $\cosh{F}$ at $x$. (An element $\beta$ of the costalk $\cosh{F}_x$ is simply a set of consistent choices $\beta_U$ of cosections over each open neighbourhood $U$ of $x$.) We topologise the display space by declaring
$$
U(\alpha) = \{ \beta \in \cosh{F}_x \ | \ x\in U, \beta_U = \alpha \}
$$
for each open $U \subset X$ and $\alpha\in \cosh{F}(U)$ to be a basis of opens. The projection
$$
p_{D\cosh{F}}: D\cosh{F} \to X
$$
with fibres the costalks is then continuous because $p_{D\cosh{F}}^{-1}U = \sum_{\alpha \in \cosh{F}(U)} U(\alpha)$ is open. It is easy to check that $D$ determines a functor $\cat{Precosh\!\downarrow\!X} \to \cat{Top \!\downarrow\! X}$.

One might imagine that $C$ and $D$ were adjoint, but this is not quite so. There is a natural map
$\epsilon_{\cosh{F}} : CD\cosh{F} \to \cosh{F}$ for any precosheaf $\cosh{F}$ given by
$$
CD\cosh{F}(U) \to \cosh{F}(U) : [\beta] \mapsto \beta_U
$$
where $[\beta]$ is the component of $p_{D\cosh{F}}^{-1}U$ containing $\beta \in \cosh{F}_x$. This is well-defined since $p_{D\cosh{F}}^{-1}U = \sum_{\alpha \in \cosh{F}(U)} U(\alpha)$ is a disjoint union of opens so each component is contained within a unique $U(\alpha)$. However, the natural map 
$$
\eta_Y : Y \to DCY : y \mapsto \{ [y] \in \pi_0(p_Y^{-1}U) \}_{U\ni x}
$$
need not be continuous. The inverse image of the basic open subset $U(\alpha)$ is the connected component $\alpha$ of $p_Y^{-1}U$ where $\alpha \in \pi_0(p^{-1}U)$. This need not be open unless $Y$ is locally-connected. 
Fortunately, the solution to this difficulty is simple. The inclusion $\cat{LCTop \!\downarrow\! X} \hookrightarrow \cat{Top \!\downarrow\! X}$ has a right adjoint $Y \mapsto \widehat{Y}$ given by the unique minimal refinement of the topology on $Y$ which is locally-connected, see \eg \cite[\S5]{funk}. Write $\widehat{D}$ for the functor $\cosh{F} \mapsto \widehat{D\cosh{F}}$.
\begin{proposition}
The functor $C:\cat{LCTop \!\downarrow\! X} \to \cat{Precosh\!\downarrow\!X}$ is left adjoint to $\widehat{D} : \cat{Precosh\!\downarrow\!X} \to \cat{LCTop \!\downarrow\! X}$.
\end{proposition}
\begin{proof}
The unit $\eta_Y : Y\to DCY$ and counit $\epsilon_{\cosh{F}} : CD\cosh{F}\to \cosh{F}$ were constructed above \emph{when $Y$ was locally-connected}. We can easily check that in this case these determine a natural isomorphism
$$
\cat{Precosh\!\downarrow\!X}\left(CY , \cosh{F} \right) \cong \cat{Top \!\downarrow\! X}\left( Y , D\cosh{F} \right).
$$
We obtain the result by composing with the natural isomorphism
$$
\cat{Top \!\downarrow\! X}\left( Y , D\cosh{F} \right) \cong \cat{LCTop \!\downarrow\! X}\left( Y , \widehat{D}\cosh{F} \right).
$$
\end{proof}

\begin{definition}
\label{spatial}
A precosheaf $\cosh{F}$ is \defn{spatial} if the basic open sets of its display space $D\cosh{F}$ are non-empty and connected.
\end{definition}

\begin{proposition}
\label{spatial display}
The following are equivalent for a precosheaf $\cosh{F}$:
\begin{enumerate}
\item $\cosh{F}$ is spatial;
\item the counit $C\widehat{D}\cosh{F}\to \cosh{F}$ is an isomorphism;
\item $\cosh{F}\cong CY$ for some locally-connected $p: Y\to X$.
\end{enumerate}
\end{proposition}
\begin{proof}
Suppose that $\cosh{F}$ is a spatial precosheaf. Then the display space $D\cosh{F}$ is locally-connected. Hence $D\cosh{F} = \widehat{D}\cosh{F}$. Moreover,
$$
CD\cosh{F}(U) = \pi_0\left(p_{D\cosh{F}}^{-1}U\right) = \pi_0\left(\sum_{\alpha\in \cosh{F}(U)} U(\alpha)\right) \cong \cosh{F}(U).
$$ 
So the counit $C\widehat{D}\cosh{F}\to \cosh{F}$ is an isomorphism. That the second statement implies the third is clear, since $\widehat{D}\cosh{F}$ is locally-connected. 

To see that the third statement implies the first, suppose, without loss of generality, that $\cosh{F} = CY$ for some locally-connected $p: Y \to X$. Let $U\subset X$ be open and $\alpha \in \cosh{F}(U)$. The corresponding basic open subset of the display space $D\cosh{F}$ is $U(\alpha) = \{\beta\in \mathcal{F}_x\ |\ x\in U , \beta_U=\alpha\}$.  Consider the unit $\eta_Y : Y \to DCY = D\cosh{F}$. It is continuous since $Y$ is locally-connected. Moreover, $\eta_Y^{-1}U(\alpha) = \alpha \in \pi_0(p^{-1}U)$ so that it, hence also $U(\alpha)$, is non-empty. Now suppose $U(\alpha) = V' + V''$ is a disjoint union of open sets. Then $\alpha = \eta_Y^{-1}U(\alpha) = \eta_Y^{-1}(V') + \eta_Y^{-1}(V'')$ is also a disjoint union of open sets. Since $\alpha$ is a connected component this is only possible if either $ \eta_Y^{-1}(V')$ or $ \eta_Y^{-1}(V'')$ is empty. It follows that either $V'$ or $V''$ is empty, and therefore that $U(\alpha)$ is connected. Hence $\cosh{F}$ is spatial.
\end{proof}
Note that, combining this with the remarks in Example \ref{component cosheaf}, a spatial precosheaf is automatically a cosheaf. Let $\cat{SpCosh\!\downarrow\!X}$ denote the full subcategory of spatial (pre)cosheaves on $X$. Another example of spatial cosheaves, required in \S \ref{stratified covers}, is provided by the next lemma.
\begin{lemma}
\label{spatial example}
Let $G : \Pi_1^{op}(X) \to \cat{Set}$ be a functor, and let $\cosh{G}$ be the precosheaf with cosections $\mathcal{G}(U) = \sum_{x\in U} G(x) \big/ \sim$ where $\sim$ is the equivalence relation generated by $a \sim a'$ if there is an op-path $\gamma$ in $U$ from $x$ to $x'$ with $a'=G(\gamma)(a)$. Then $\cosh{G}$ is spatial, in particular it is actually a cosheaf.
\end{lemma}
\begin{proof}
Using Lemma \ref{connectivity lemma} we find that the costalk $\cosh{G}_x = G(x)$. Therefore the display space $D\cosh{G}$ is $\sum_{x \in X} G(x)$ with the topology generated by the basic open sets $U(\alpha)$ for $\alpha \in \cosh{G}(U)$ where $U\subset X$ is open. Note that the cosection $\alpha \in \cosh{G}(U)$ is represented (non-uniquely) by some $a\in G(y)$ where $y\in U$. In terms of this representative, $U(\alpha)  = \{ b \in G(x) \ | \ x\in U, b \sim a \}$. It follows that the basic open sets are non-empty and connected. Hence $\cosh{G}$ is a spatial precosheaf.
\end{proof}

We have already seen that the display space $D\cosh{F}$ of a spatial cosheaf is locally-connected, it also has other good properties. The idea of a \defn{complete spread} was introduced by Fox in \cite{fox} to give a purely topological notion of a branched cover. A map $p_Y:Y \to X$ is a \defn{spread} if the set of connected components of $p_Y^{-1}U$ for open $U$ in $X$ forms a basis for the topology of $Y$. It is a \defn{complete spread} if whenever we make a consistent choice of component $\alpha_U \subset p_Y^{-1}U$ for each neighbourhood $U$ of some fixed $x\in X$ --- consistent meaning that $\alpha_U \subset \alpha_V$ whenever $U \subset V$ --- then the intersection $\cap_{U \ni x} \alpha_U \neq \emptyset$. If, in addition, there is a unique point in this intersection we say $Y$ is \defn{uniquely-complete}. More succinctly, a spread $Y$ is uniquely-complete if and only if
$$
p_Y^{-1}x \to \lim_{U \ni x} \pi_0\left(p_Y^{-1}U \right): y \mapsto \left\{ [y] \in \pi_0(p_Y^{-1}U) \right\}_{U\ni x}
$$
is a bijection for each $x\in X$. Write $\cat{UCS\!\downarrow\!X}$ for the category of locally-connected, uniquely-complete spreads over $X$. The maps are continuous maps over $X$. 

\begin{corollary}
If $\cosh{F}$ is a spatial cosheaf then $D\cosh{F}$ is a locally-connected, uniquely-complete spread.
\end{corollary}
\begin{proof}
Assume $\cosh{F}$ is a spatial cosheaf on $X$. Let $U\subset X$ be open. Then
$p_{D\cosh{F}}^{-1}U = \sum_{\alpha \in \cosh{F}(U)} U(\alpha)$. The right hand side is the decomposition into connected components. By definition these form a basis of the topology of $D\cosh{F}$, which is thus a spread. As remarked above, $D\cosh{F}$ is locally-connected, and 
$$
\lim_{U\ni x} \pi_0\left(p_{D\cosh{F}}^{-1}U \right) =
\lim_{U\ni x} \pi_0\left(\sum_{\alpha \in \cosh{F}(U)} U(\alpha) \right)
\cong \lim_{U\ni x} \cosh{F}(U) = \cosh{F}_x
$$
so that $D\cosh{F}$ is uniquely-complete. 
\end{proof}

\begin{proposition}
The unit $Y \to \widehat{D}CY$ is a homeomorphism if and only if  $Y$ is a locally-connected, uniquely-complete spread over $X$. 
\end{proposition}
\begin{proof}
 If $Y$ is locally-connected then by  Proposition \ref{spatial display} $CY$ is a spatial cosheaf.  Hence $DCY$ is locally-connected. Hence $\widehat{D}CY = DCY$. Consider the map $Y \to DCY$. The restriction to the fibre over $x$ is
$$
p_Y^{-1}x \to  \lim_{U \ni x} \pi_0\left(p_Y^{-1}U \right) = CY_x = p_{DCY}^{-1}x
$$
which is a bijection precisely when $Y$ is uniquely-complete. When this holds, the map is a homeomorphism if and only if the inverse images of the basis $\{U(\alpha)\}$  for the topology of $DCY$ are a basis for the topology of $Y$. The inverse image of $U(\alpha)$, where $\alpha \in CY(U) = \pi_0(p_Y^{-1}U)$, is simply the connected component $\alpha$ of $p_Y^{-1}U$. So the map is a homeomorphism precisely when these components form a basis, \ie precisely when $Y$ is a spread.
\end{proof}
We have shown that the adjoint functors $C$ and $\widehat{D}$ restrict to an equivalence between the full subcategories of locally-connected, uniquely-complete spreads and of spatial cosheaves. The situation is summarised in the following commutative diagram (the vertical arrows are the inclusions).
\[
\xymatrix{
\cat{LCTop \!\downarrow\! X} \ar[dr]^{C} & \cat{Precosh\!\downarrow\!X} \ar[l]_{\widehat{D}} \\
\cat{UCS\!\downarrow\!X} \ar[u] \ar@{<->}[r]^\sim
& \cat{SpCosh\!\downarrow\!X} \ar[u] 
}
\]
\begin{remark}
The composite $C\widehat{D} : \cat{Precosh\!\downarrow\!X} \to \cat{SpCosh\!\downarrow\!X}$ is a `spatial cosheafification' functor.  That is, for any precosheaf $\cosh{F}$ there is a map $C\widehat{D}\cosh{F} \to \cosh{F}$ with the universal property that for any spatial cosheaf $\cosh{E}$ and map $\varphi:\cosh{E}\to \cosh{F}$  of precosheaves there is a unique factorisation:
$$
\xymatrix{
& C\widehat{D}\cosh{F} \ar[d] \\
\cosh{E} \ar[r]_\varphi \ar@{-->}[ur]^{\exists !} & \cosh{F}
}
$$
\end{remark}

We will not go into a full discussion of the functoriality of cosheaves here. However note that a cosheaf can be restricted to a subspace $\imath: Y \to X$ by defining
$$
\imath^*\cosh{F}(V) = \lim_{U \supset V} \cosh{F}(U).
$$
This corresponds to restricting the corresponding locally-connected uniquely-complete spread $D\cosh{F}$ to $Y$. The restriction of a cosheaf to a point is simply the costalk.

\bibliographystyle{alpha}
\bibliography{pohomotopy}

 \end{document}